\documentclass[11pt]{article}
\usepackage{amsmath}
\usepackage{amsfonts}
\usepackage{amssymb}
\usepackage{hyperref}
\usepackage{indentfirst}
\usepackage{mathrsfs}

\usepackage[top=1in,bottom=1in,left=1.25in,right=1.25in]{geometry}
\def\<{\langle}
\def\>{\rangle}

\date{}
\allowdisplaybreaks

\begin{document}

\renewcommand{\baselinestretch}{1.2}
\renewcommand{\arraystretch}{1.0}

\title{\bf Generalized Yetter-Drinfel'd module categories for regular multiplier Hopf algebras}
\author
{
  \textbf{Tao Yang} \footnote{Corresponding author. College of Science, Nanjing Agricultural University, Nanjing 210095, Jiangsu, CHINA.
             E-mail: tao.yang.seu@gmail.com}, \,
  \textbf{Xuan Zhou} \footnote{Department of Mathematics, Jiangsu Institute of Education, Nanjing 210013, Jiangsu, CHINA.
             E-mail: zhouxuanseu@126.com}
}

\maketitle

\begin{center}
\begin{minipage}{12.cm}

 {\bf Abstract }
 For a commutative regular multiplier Hopf algebra $A$, the Yetter-Drinfel'd module category ${}_{A}\mathcal{YD}^{A}$
 is equivalent to the centre $Z({}_{A}\mathcal{M})$ of the unital left $A$-module category ${}_{A}\mathcal{M}$.
 Then we introduce the generalized $(\alpha, \beta)$-Yetter-Drinfel'd module categories
 ${}_{A}\mathcal{GYD}^{A}(\alpha, \beta)$, which are treated as components of a braided $T$-category.
 Especially when $A$ is a coFrobenius Hopf algebra,
 ${}_{A}\mathcal{YD}^{A}(\alpha, \beta)$ is isomorphic to the unital $\widehat{A} \bowtie A(\alpha, \beta)$-module category
 ${}_{\widehat{A} \bowtie A(\alpha, \beta)}\mathcal{M}$.
 Finally for a Yetter-Drinfel'd $A$-module algebra $H$, we introduce Yetter-Drinfel'd $(H, A)$-module category, which is a monoidal.
\\

\\


\end{minipage}
\end{center}
\normalsize

\section{Introduction}
\def\theequation{\thesection.\arabic{equation}}
\setcounter{equation}{0}

 A Yetter-Drinfel'd module (or crossed bimodule, Yang-Baxter module) over Hopf algebra is a module and a comodule
 satisfying a certain compatibility condition.
 The main feature of this definition is that Yetter-Drinfel'd modules form a pre-braided monoidal category.
 Under favourable conditions (e.g. the antipode of Hopf algebra is bijective), the category is even braided (or quasisymmetric).
 Via (pre-)braiding structure, the notion of Yetter-Drinfel'd module plays a part in the relations between quantum group and knot theory.

 As Hopf algebra is finite dimensional, the Yetter-Drinfel'd module category has some special properties,
 e.g., it is isomorphic to the category of left modules over the Drinfel'd double of this Hopf algebra.
 This property no longer holds in the infinite dimensional case.
 However by the multiplier Hopf algebra theory, the isomorphism exists when the Hopf algebra is coFrobenius (see \cite{Zh99}).
 Now, multiplier Hopf algebra becomes a valuable tool to deal with infinite dimensional Hopf algebra cases.

 In \cite{YW11}, new kinds of Yetter-Drinfeld modules over a regular multiplier Hopf algebra were introduced via extended module structure.
 We also use them to construct new braided crossed categories, which are related to homotopy invariants.

 In this paper, we continue to consider the (left-right) Yetter-Drinfel'd module category $_{A}\mathcal {YD}^{A}$  defined in \cite{YW11},
 and show that it is isomorphic to the center of unital left $A$-module category.
 Then we consider two kinds of generalizations of Yetter-Drinfel'd module categories and get some properties.

 The paper is organized in the following way.
 In section 2, we recall some notions which we will use in the following, such as
 multiplier Hopf algebras, extended modules, comodules and centre of a strict tensor category.

 In section 3, we mianly show the relationship between the Yetter-drinfel'd module category ${}_{A}\mathcal{YD}^{A}$
 and the centre $Z({}_{A}\mathcal{M})$ of the unital $A$-module category ${}_{A}\mathcal{M}$.

 In section 4, we introduced the generalized $(\alpha, \beta)$-Yetter-Drinfel'd module categories
 ${}_{A}\mathcal{GYD}^{A}(\alpha, \beta)$, which are treated as components of a braided $T$-category.
 Especially when $A$ is a coFrobenius Hopf algebra,
 ${}_{A}\mathcal{YD}^{A}(\alpha, \beta)$ is isomorphic to the unital $\widehat{A} \bowtie A(\alpha, \beta)$-module category
 ${}_{\widehat{A} \bowtie A(\alpha, \beta)}\mathcal{M}$.

 In section 5, let $H$ be a Yetter-Drinfel'd $A$-module algebra, we introduce Yetter-Drinfel'd $(H, A)$-module category, which is a monoidal.

\section{Preliminaries}
\def\theequation{\thesection.\arabic{equation}}
\setcounter{equation}{0}

 We begin this section with a short introduction to multiplier Hopf algebras.

 Throughout this paper, all spaces we considered are over a fixed field $K$ (such as field $\mathbb{C}$ of complex numbers).
 Algebras may or may not have units, but always should be non-degenerate,
 i.e., the multiplication maps (viewed as bilinear forms) are non-degenerate.

 For an algebra $A$, the multiplier algebra $M(A)$ of $A$ is defined as the largest algebra with unit
 in which $A$ is a dense ideal (see the appendix in \cite{V94}).

 Now, we recall the definitions of a multiplier Hopf algebra (see \cite{V94} for details).
 A comultiplication on an algebra $A$ is a homomorphism $\Delta: A \longrightarrow M(A \otimes A)$
 such that $\Delta(a)(1 \otimes b)$ and $(a \otimes 1)\Delta(b)$ belong to $A\otimes A$ for all $a, b \in A$.
 We require $\Delta$ to be coassociative in the sense that
 \begin{eqnarray*}
 (a\otimes 1\otimes 1)(\Delta \otimes \iota)(\Delta(b)(1\otimes c))
 = (\iota \otimes \Delta)((a \otimes 1)\Delta(b))(1\otimes 1\otimes c)
 \end{eqnarray*}
 for all $a, b, c \in A$ (where $\iota$ denotes the identity map).

 A pair $(A, \Delta)$ of an algebra $A$ with non-degenerate product and a comultiplication $\Delta$ on $A$ is called a multiplier Hopf algebra,
 if the linear maps $T_{1}, T_{2}: A\otimes A \longrightarrow A\otimes A$ defined by
 \begin{eqnarray}
 T_{1}(a\otimes b)=\Delta(a)(1 \otimes b), \qquad  T_{2}(a\otimes b)=(a \otimes 1)\Delta(b)
 \end{eqnarray}
 are bijective.

 A multiplier Hopf algebra $(A, \Delta)$ is called regular if $(A, \Delta^{cop})$ is also a multiplier Hopf algebra,
 where $\Delta^{cop}$ denotes the co-opposite comultiplication defined as $\Delta^{cop}=\tau \circ \Delta$ with $\tau$ the usual flip map
 from $A\otimes A$ to itself (and extended to $M(A\otimes A)$). In this case,
 $\Delta(a)(b \otimes 1) \mbox{ and } (1 \otimes a)\Delta(b) \in A \otimes A$
 for all $a, b\in A$.

 Multiplier Hopf algebra $(A, \Delta)$ is regular if and only if the antipode $S$ is bijective from $A$ to $A$
 (see \cite{V98}, Proposition 2.9). In this situation, the comultiplication is also determined by the bijective maps
 $T_{3}, T_{4}: A\otimes A \longrightarrow A\otimes A$ defined as follows
 \begin{eqnarray}
 && T_{3}(a \otimes b)=\Delta(a)(b \otimes 1), \qquad T_{4}(a\otimes b)=(1 \otimes a)\Delta(b)
 \end{eqnarray}
 for all $a, b\in A$.

 We will use the adapted Sweedler notation for regular multiplier Hopf algebras (see \cite{V08}).
 We will e.g., write $\sum a_{(1)} \otimes a_{(2)}b$ for $\Delta(a)(1 \otimes b)$
 and $\sum ab_{(1)} \otimes b_{(2)}$ for $(a \otimes 1)\Delta(b)$, sometimes we omit the $\sum$.

 Define two linear operators $\mathcal{T}$ and $\mathcal{T'}$ acting on $A \otimes A$ introduced by the formulae
 \begin{eqnarray*}
 && \mathcal{T}(a \otimes b) = b_{(2)} \otimes aS(b_{(1)})b_{(3)}, \\
 && \mathcal{T'}(a \otimes b) = b_{(1)} \otimes S(b_{(2)})ab_{(3)}
 \end{eqnarray*}
 for any $a, b\in A$.

 This two operators above are well-defined, since
 \begin{eqnarray*}
 && \mathcal{T}(a \otimes b) = T_{4}(S\otimes \iota)T_{3}(\iota\otimes S^{-1})\tau (a\otimes b), \\
 && \mathcal{T'}(a \otimes b) = (\iota \otimes S) T_{4} \tau (\iota\otimes S^{-1}) T_{4} (a\otimes b).
 \end{eqnarray*}
 They are obviously invertible, and the inverses can be written as follows
 \begin{eqnarray*}
 && \mathcal{T}^{-1}(a \otimes b) = b S^{-1}(a_{(3)}) a_{(1)} \otimes a_{(2)}, \\
 && \mathcal{T'}^{-1}(a \otimes b) = a_{(3)} b S^{-1}(a_{(2)}) \otimes a_{(1)}.
 \end{eqnarray*}

 For any $a, b \in A$, $\mathcal{T} \circ T_{2} = T_{4}$. If $A$ is commutative, then $\mathcal{T'}=\tau$,
 and if $A$ is cocommutative, then $\mathcal{T}=\tau$.

 \textbf{Proposition \thesection.1}
 Operators $\mathcal{T}$ and $\mathcal{T'}$ satisfy the braided equation
 \begin{eqnarray*}
 && (\mathcal{T} \otimes \iota)(\iota \otimes \mathcal{T})(\mathcal{T} \otimes \iota)
  =(\iota \otimes \mathcal{T})(\mathcal{T} \otimes \iota)(\iota \otimes \mathcal{T}), \\
 && (\mathcal{T'} \otimes \iota)(\iota \otimes \mathcal{T'})(\mathcal{T'} \otimes \iota)
  =(\iota \otimes \mathcal{T'})(\mathcal{T'} \otimes \iota)(\iota \otimes \mathcal{T'}).
 \end{eqnarray*}

 \subsection{Extended modules}

 Let $A$ be an (associative) algebra. Suppose $X$ is a left $A$-module
 with the module structure map $\cdot : A\otimes X\longrightarrow  X$. We will always assume that the module is non-degenerate,
 this means that $x=0$ if $x\in X$ and $a\cdot x=0$ for all $a\in A$.
 If the module is unital (i.e., $A\cdot X=X$), then we can get an extension of the module structure to $M(A)$,
 this means that we can define $f\cdot x$, where $f\in M(A)$ and $x\in X$.
 In fact, since $x\in X=A\cdot X$, then $x= \sum_{i} a_{i}\cdot x_{i}$ and $f\cdot x=\sum_{i} (fa_{i})\cdot x_{i}$.
 In this setting, we can easily get $1\cdot x=x$, where $1$ is the unit of multiplier algebra $M(A)$.
\smallskip

 Denote by $Y$ the space of linear maps $\rho: A\rightarrow X$ satisfying $\rho(aa')=a\cdot \rho(a')$ for all $a, a'\in A$.
 Then $Y$ becomes a left $A$-module if we define $a\cdot \rho$ for $a\in A$ and $\rho \in Y$ by $(a\cdot \rho)(a')=\rho(a'a)$.
 For any $x\in X$, we have an element $\rho_{x} \in Y$ defined by $\rho_{x}(a)=a\cdot x$ when $a\in A$.
 The injective map $x \mapsto \rho_{x}$ is a map of left modules, from this $X$ becomes a submodule of $Y$.
 We call $Y$ the extended module of $X$.

 From the definition, we have $A\cdot Y\subseteq X$, and $A\cdot Y=X$ if $A\cdot X=X$.
 Since $A^{2}=A$, $Y$ is still non-degenerate. If $A$ has a unit, then $Y=X$, in the other case,
 mostly $Y$ is strictly bigger than $X$.

 Because the module $Y$ is the completion of the original left $A$-module $X$ for the strict topology (see \cite{V08}),
 sometimes it is also called the complete module (or completion) of $X$, and denoted $Y$ as $M_{l}(X)$.

 We can also do the same for a right $A$-module $X$ and denote the extended module of the right $A$-module $X$ as $M_{r}(X)$.

 Let $X$ be a non-degenerate $A$-bimodule.
 Denote by $Z$ the space of pair $(\lambda, \rho)$ of linear maps from $A$ to $X$  satisfying
 $a\cdot \lambda(a')=\rho(a)\cdot a'$ for all $a, a' \in A$.
 From the non-degeneracy, it follows that
 $\rho(aa')=a\cdot \rho(a')$ and $\lambda(aa')=\lambda(a)\cdot a'$ for all $a, a' \in A$.
 Also $\rho$ is completely determined by $\lambda$ and vice versa.
 We can consider $Z$ as the intersection of two extensions of $X$ (as a left and a right modules).

 $Z$ becomes an $A$-bimodule, if we define $a\cdot z$ and $z\cdot a$ for $a\in A$ and $z=(\lambda, \rho) \in Z$ by
 $a\cdot z = (a\cdot \lambda(\cdot), \rho(\cdot a))$ and $z\cdot a = (\lambda(a \cdot), \rho(\cdot)\cdot a)$.
 If we define $(\lambda_{x}, \rho_{x})$ for $x\in X$ by $\lambda_{x}(a)=x\cdot a$ and $\rho_{x}(a)=a\cdot x$,
 we get $X$ as a submodule of $Z$. We call $Z$ the extended module of $A$-bimodule $X$, and denote it as $M_{0}(X)$.

 Let $V$ be a vector space and $X = V\otimes A$. We consider the left and right module actions of $A$ on $X$ respectively given by
 \begin{eqnarray*}
  a \cdot (v \otimes a')= v\otimes aa' \quad \mbox{and} \quad (v \otimes a')\cdot a= v\otimes a'a
 \end{eqnarray*}
 for $a, a'\in A$ and $v\in V$.
 We denote the extended module of left $A$-module $X$ as $M_{l}(V \otimes A)$, and extended module of the right one as $M_{r}(V \otimes A)$.
 This two actions are compatible, and make $X$ an $A$-bimodule, we denote its extend module as $M_{0}(V \otimes A)$.

 As vector spaces, $M_{l}(V \otimes A)$ and $M_{0}(V \otimes A)$ are bigger than $V \otimes A$,
 since $V \otimes R(A) \subseteq M_{l}(V \otimes A)$ and $V \otimes M(A) \subseteq M_{0}(V \otimes A)$,
 where $R(A)$ is the right multipliers of $A$ and $M(A)$ is the multipliers of $A$ (see Appendix in \cite{V94}).
 If $V$ is finite-dimensional, then $M_{l}(V \otimes A) = V \otimes R(A)$.

 \subsection{Comodules of a regular multiplier Hopf algebra}

 Let $A$ be a regular multiplier Hopf algebra and $V$ a vector space.
 A right coaction of $A$ on $V$ is an injective linear map $\Gamma: V \rightarrow M_{r}(V \otimes A)$,
 where $M_{r}(V \otimes A)$ is the extended module of the right $A$-module $V\otimes A$, such that
 $$(\Gamma \otimes \iota)\Gamma=(\iota \otimes \Delta)\Gamma.$$

 We call $V$ a right $A$-comodule.
 Of course, we need to give a precise meaning to the last equation. Indeed, similar to Proposition 2.10 in \cite{V08},
 the map $\iota \otimes \Delta$ and $\Gamma \otimes \iota$, defined on $V\otimes A$, have natural extensions to maps
 from $M_{r}(V \otimes A)$ to $M_{r}(V \otimes A \otimes A)$.

 For any $v \in V$, $\Gamma(v)\in M_{r}(V \otimes A)$ and $\Gamma(v)(1\otimes a) \in V \otimes A$.
 We will also use the adapted version of the Sweedler notation for this coaction, and write
 $\Gamma(v)(1\otimes a)=v_{(0)}\otimes v_{(1)}a$ for all $a\in A$ and $v\in V$.

 Similarly, we can also define comodules by the other two kind of extended modules.
 These definitions are essentially the same, except for the ranges of the coactions.

 \subsection{Centre of a strict tensor category}

 Let $(\mathcal{C}, \otimes, I)$ be a strict tensor category.
 We first recall the definition of the centre $Z(\mathcal{C})$ and some of its properties (see \cite{K}).

 An object of $Z(\mathcal{C})$ is a pair $(V, C_{-, V})$, where $V$ is an object of $\mathcal{C}$ (i.e., $V\in \mathcal{C}$)
 and $C_{-, V}$ is a family of natural isomorphisms $C_{-, V}: X\otimes V\rightarrow V\otimes X$
 defined for all objects $X$ in $\mathcal{C}$ such that for all objects $X, Y$ in $\mathcal{C}$ we have
 $C_{X\otimes Y, V} = (C_{X, V} \otimes \iota)(\iota \otimes C_{Y, V})$.
 The naturality means that for any morphism $f: X \rightarrow Y$ in $\mathcal{C}$,
 $(\iota \otimes f)C_{X, V} = C_{Y, V} (f \otimes \iota)$.

 A morphism from $(V, C_{-, V})$ to $(W, C_{-, W})$ is a morphism $f: V\rightarrow W$ such that for
 each object $X$ of $\mathcal{C}$ we have $(f\otimes \iota)C_{X, V} = C_{X, W}(\iota \otimes f)$.

 Let $(\mathcal{C}, \otimes, I)$ be a strict tensor category. then $Z(\mathcal{C})$ is a strict tensor category,
 where
 \begin{enumerate}
 \item[(1)] the unit is $(I, \iota)$;

 \item[(2)] the tensor product of $(V, C_{-, V})$ and $(W, C_{-, W})$ is given by
 $(V, C_{-, V}) \otimes (W, C_{-, W}) = (V\otimes W, C_{-, V\otimes W})$,
 where $C_{-, V\otimes W}: X\otimes V\otimes W \rightarrow V\otimes W\otimes X$
 is the morphism of $\mathcal{C}$ defined for all objects $X$ in $\mathcal{C}$ by
 $C_{X, V\otimes W} = (\iota\otimes C_{X, W})(C_{X, V}\otimes \iota)$;

 \item[(3)] and the braiding is given by $C_{V, W}: (V, C_{-, V}) \otimes (W, C_{-, W}) \rightarrow (W, C_{-, W}) \otimes (V, C_{-, V})$.
 \end{enumerate}

\section{A categorical interpretation of Yetter-Drinfel'd modules}
\def\theequation{\thesection.\arabic{equation}}
\setcounter{equation}{0}

 Let $(A, \Delta, \varepsilon, S)$ be a regular multiplier Hopf algebra.
 We first recall from \cite{YW11} that a (left-right) Yetter-Drinfel'd module over $A$ is a vector space $V$
 satisfying the following conditions
 \begin{itemize}
 \item $(V, \cdot)$ is an unital left $A$-module, i.e., $A\cdot V = V$.

 \item $(V, \Gamma)$ is a (right) $A$-comodule, where
 $\Gamma: V\rightarrow M_{0}(V\otimes A)$ denotes the right coaction of $A$ on $V$,
 $M_{0} (V\otimes A)$ denote the extended module shown in section 2.

 \item $\Gamma$ and $\cdot$ satisfy the following compatible conditions
 \begin{eqnarray}
 && (a\cdot v)_{(0)} \otimes (a\cdot v)_{(1)} a' = a_{(2)}\cdot v_{(0)} \otimes a_{(3)} v_{(1)} S^{-1}(a_{(1)}) a'. \label{2}
 \end{eqnarray}
 \end{itemize}

 The Yetter-Drinfel'd module category $_{A}\mathcal {YD}^{A}$ is defined as follows.
 The objects in $_{A}\mathcal {YD}^{A}$ are left-right Yetter-Drinfel'd modules and the morphisms are linear maps
 that intertwine with the left action and the right coaction of $A$ on $M$,
 i.e., the morphisms between two objects are left $A$-linear and right $A$-colinear maps.
 More precisely, let $V, W \in {}_{A}\mathcal {YD}^{A}$ and $f: V\rightarrow W$ be a morphism,
 then
 \begin{eqnarray}
 && f(a\cdot v) = a\cdot f(v), \nonumber \\
 && \Gamma_{W}\circ f(v)(1\otimes a') = (f\otimes \iota)\Gamma_{V}(v)(1\otimes a'),
 \end{eqnarray}
 for all $a, a' \in A$ and $v \in V$, where $\Gamma_{W}$ ($\Gamma_{V}$) is the right coaction on $W$ ($V$).

 $_{A}\mathcal{YD}^{A}$ is a braided monoidal category in the following way. For $V, W \in {}_{A}\mathcal {YD}^{A}$,
 the action and coaction of $A$ on $V\otimes W$ is given by
 \begin{eqnarray}
 && a\cdot (v\otimes w)=a_{(1)}\cdot v \otimes a_{(2)}\cdot w, \\
 && \Gamma(v\otimes w)(1\otimes a')= v_{(0)} \otimes w_{(0)} \otimes w_{(1)} v_{(1)}a'.
 \end{eqnarray}
 for any $v\in V, w \in W$ and $a, a' \in A$.
 \\

 Let $A$ be a regular multiplier Hopf algebra and denote by ${}_{A}\mathcal{M}$ the unital left $A$-module category.
 The module structure is same as what in Yetter-Drinfel'd category.
 We will try to show the equivalence between $Z({}_{A}\mathcal{M})$ and ${}_{A}\mathcal{YD}^{A}$.
 \\

 \textbf{Proposition \thesection.1}  Let $V$ be an object in ${}_{A}\mathcal{YD}^{A}$.
 For any left $A$-module $X$, define $C_{X, V}(x\otimes v)=v_{(0)} \otimes v_{(1)}\cdot x$ where $x\in X$ and $v\in V$,
 then $(V, C_{-, V}) \in Z({}_{A}\mathcal{M})$.

 \emph{Proof.}
 Indeed, $C_{X, V}$ is well-defined, and clearly it is invertible
 $C^{-1}_{X, V}(v\otimes x) = S(v_{(1)})\cdot x \otimes v_{(0)}$.
 For any unital left $A$-module $Y$, $x\in X$, $y\in Y$, and $v\in V$,
 \begin{eqnarray*}
 C_{X\otimes Y, V}(x\otimes y\otimes v)
 &=& v_{(0)} \otimes v_{(1)}\cdot x \otimes v_{(2)}\cdot y \\
 &=& C_{X, V}(x\otimes v_{(0)}) \otimes v_{(1)}\cdot y \\
 &=& (C_{X, V}\otimes \iota)(x\otimes v_{(0)} \otimes v_{(1)}\cdot y) \\
 &=& (C_{X, V}\otimes \iota)(\iota \otimes C_{Y, V})(x\otimes y \otimes v),
 \end{eqnarray*}
 we get that $C_{X\otimes Y, V} = (C_{X, V}\otimes \iota)(\iota \otimes C_{Y, V})$.
 Let $f: X\rightarrow Y$ be a left $A$-module map,
 \begin{eqnarray*}
 C_{Y, V}(f\otimes \iota)(x\otimes v)
 &=& C_{Y, V}(f(x)\otimes v)
 = v_{(0)} \otimes v_{(1)}\cdot f(x) \\
 &=& v_{(0)} \otimes f(v_{(1)}\cdot x)
 = (\iota \otimes f)C_{X, V}(x\otimes v),
 \end{eqnarray*}
 thus $C_{-, V}$ is natural. $C_{X, V}$ is a left $A$-module map, since
 \begin{eqnarray*}
 C_{X, V}(a\cdot (x\otimes v))
 &=& C_{X, V}(a_{(1)}\cdot x\otimes a_{(2)}\cdot v)
 = (a_{(2)}\cdot v)_{(0)} \otimes (a_{(2)}\cdot v)_{(1)} a_{(1)} x \\
 &=& a_{(1)}\cdot v_{(0)}\otimes a_{(2)} v_{(1)} x
 = a\cdot C_{X, V}(x\otimes v),
 \end{eqnarray*}
 where the third equation holds because of the Yetter-Drinfel'd compatible condition.
 $\hfill \blacksquare$
 \\

 Let $f: V\rightarrow W$ be a morphism in ${}_{A}\mathcal{YD}^{A}$, i.e., $f$ is a module and comodule map.
 Then
 \begin{eqnarray*}
 (f\otimes \iota)C_{X, V}(x\otimes v)
 &=& f(v_{(0)}) \otimes v_{(1)}\cdot x
 = f(v)_{(0)} \otimes f(v)_{(1)}\cdot x \\
 &=& C_{X, W}(x\otimes f(v))
 = C_{X, W}(\iota\otimes f)(x\otimes v),
 \end{eqnarray*}
 so $(f\otimes \iota)C_{X, V} = C_{X, W}(\iota\otimes f)$.

 We can define a functor $G: {}_{A}\mathcal{YD}^{A} \rightarrow Z({}_{A}\mathcal{M})$ as follows,
 \begin{eqnarray*}
 G(V)=(V, C_{-, V}), \qquad G(f)=f,
 \end{eqnarray*}
 where $C_{X, V}(x\otimes v)=v_{(0)} \otimes v_{(1)}\cdot x$, $f: V\rightarrow W$ is a morphism in ${}_{A}\mathcal{YD}^{A}$.
 Also we can get that $G$ preserves tensor product, since for $V, W \in {}_{A}\mathcal{YD}^{A}$, $X\in {}_{A}M$,
 $x\in X, v\in V, w\in W$,
 \begin{eqnarray*}
 C_{X, V\otimes W}(x\otimes v\otimes w)
 &=& v_{(0)} \otimes w_{(0)} \otimes w_{(1)}v_{(1)}\cdot x \\
 &=& v_{(0)} \otimes C_{X, W}(v_{(1)}\cdot x \otimes w) \\
 &=& (\iota\otimes C_{X, W})(v_{(0)} \otimes v_{(1)}\cdot x \otimes w) \\
 &=& (\iota\otimes C_{X, W})(C_{X, V} \otimes \iota)(x \otimes v \otimes w).
 \end{eqnarray*}

 \textbf{Proposition \thesection.2} Let $A$ be a regular multiplier Hopf algebra.
 For any object $(V, C_{-, V})$ in $Z({}_{A}\mathcal{M})$, if $M_{r}(V\otimes A) = M_{0}(V\otimes A)$ as vector spaces,
 then $V$ is a Yetter-Drinfel'd module over $A$.

 \emph{Proof.}
 For any object $(V, C_{-, V})$ in $Z({}_{A}\mathcal{M})$, we can endow $V$ with the structure of a Yetter-Drinfel'd module.
 For any unital left $A$-module $X$ and any element $x\in X$, there exist a local unit $e\in A$ such that $e\cdot x =x$,
 define the left $A$-module map $\bar{x}: A\rightarrow X$ by $\bar{x}(a)=a\cdot x$.
 For each $v\in V$ and $a\in A$, by the natruality of $C_{X, V}$, we have
 \begin{eqnarray*}
 && C_{X, V}(x\otimes v) = C_{X, V}(e\cdot x\otimes v)\\
 &=& C_{X, V}(\bar{x} \otimes \iota)(e\otimes v)
 = (\iota \otimes \bar{x})C_{A, V}(e\otimes v)
 \in V \otimes A\cdot x.
 \end{eqnarray*}
 Set $X=A$ with the multiplication as module action,
 then $C_{A, V}(ba\otimes v) = C_{A, V}(b\otimes v)(1\otimes a)$ for all $a, b\in A$ and $v\in V$.

 If we treat $A\otimes V$ and $V\otimes A$ as right $A$-module
 with structures acting by multiplication on the $A$-component,
 then obviously the right module action is compatible with the left one, and $C_{A, V}(b\otimes v)\cdot a = C_{A, V}((b\otimes v)\cdot a)$.
 There is a natural extension $M_{r}(A\otimes V)\longrightarrow M_{r}(V\otimes A)$, denoted by $C'_{A, V}$,
 $C'_{A, V}(y)\cdot a = C_{A, V}(y\cdot a), y\in M_{r}(A\otimes V)$ for all $a\in A$.
 So by the assumption we can get $C'_{A, V}(1\otimes v) \in M_{0}(V\otimes A)$, where $1$ is the unit of $M(A)$.

 If we define $\rho: V\rightarrow M_{0}(V\otimes A)$ by $\rho(v) = C'_{A, V}(1\otimes v)$,
 and write $C_{A, V}(a\otimes v) = C'_{A, V}(1\otimes v) \cdot a = \rho(v)(1\otimes a) = v_{(0)} \otimes v_{(1)} a$,
 then $V$ is a right $A$-comodule, since by the definition of $Z(\mathcal{C})$,
 \begin{eqnarray*}
 && C_{A\otimes A, V}(x\otimes y\otimes v)
 = v_{(0)} \otimes v_{(1)}\cdot (x\otimes y)
 = v_{(0)} \otimes v_{(1)(1)}\cdot x\otimes v_{(1)(2)}\cdot y \\
 && \qquad \quad = (C_{A, V}\otimes \iota)(\iota \otimes C_{A, V})(x\otimes y\otimes v)
 = v_{(0)(0)} \otimes v_{(0)(1)}\cdot x \otimes v_{(1)}\cdot y,
 \end{eqnarray*}
 then $v_{(0)} \otimes v_{(1)(1)} x\otimes v_{(1)(2)} y = v_{(0)(0)} \otimes v_{(0)(1)} x \otimes v_{(1)} y$
 for all $x, y\in A$, which proves that $V$ is a right $A$-comodule.

 $V$ is also a Yetter-Drinfel'd module. Since $C_{X, V}$ is a natural isomorphism in ${}_{A}\mathcal{M}$,
 $C_{X, V}$ is left $A$-linear, i.e., for all $a\in A$, $x\in X$ and $v\in V$,
 $a\cdot C_{X, V}(x\otimes v)=C_{X, V}(a\cdot (x\otimes v))$.
 \begin{eqnarray*}
 && a\cdot C_{X, V}(x\otimes v)
 = a\cdot (v_{(0)} \otimes v_{(1)}\cdot x)
 = a_{(1)}\cdot v_{(0)} \otimes a_{(2)}v_{(1)}\cdot x, \\
 && C_{X, V}(a\cdot (x\otimes v))
 = C_{X, V}(a_{(1)}\cdot x\otimes a_{(2)}\cdot v)
 = (a_{(2)}\cdot v)_{(0)} \otimes (a_{(2)}\cdot v)_{(1)} a_{(1)}\cdot x.
 \end{eqnarray*}
 From above, we can get the compatible condition of the Yetter-Drinfel'd module, i.e., for any $a'\in A$,
 $
 (a_{(2)}\cdot v)_{(0)} \otimes (a_{(2)}\cdot v)_{(1)} a_{(1)} a'
 =a_{(1)}\cdot v_{(0)} \otimes a_{(2)}v_{(1)} a'.
 $
 $\hfill \blacksquare$
 \\

 Assume that $f: (V, C_{-, V}) \rightarrow (W, C_{-, W})$ is a morphism in $Z({}_{A}\mathcal{M})$.
 Of course $f$ is a module map. Since $(f\otimes \iota)C_{A, V} = C_{A, W}(\iota\otimes f)$, $f$ is also a comodule map.
 Hence, we can define a functor $F$ from $Z({}_{A}\mathcal{M})$ to ${}_{A}\mathcal{YD}^{A}$ as follows:
 \begin{eqnarray*}
 F((V, C_{-, V}))=V, \qquad F(f)=f.
 \end{eqnarray*}
 We can easily check that $F$ preserve the tensor product. For any $v\in V$ and $w\in W$,
 \begin{eqnarray*}
 \rho(v\otimes w)(1\otimes a)
 &=& C_{A, V\otimes W}(a\otimes v\otimes w)
 = (\iota \otimes C_{A, W})(C_{A, V} \otimes \iota)(a\otimes v\otimes w) \\
 &=& (\iota \otimes C_{A, W})(v_{(0)} \otimes v_{(1)} a \otimes w)
 = v_{(0)} \otimes w_{(0)} \otimes w_{(1)}v_{(1)} a,
 \end{eqnarray*}
 so $F((V\otimes W, C_{-, V\otimes W})) = V\otimes W$.
 \\

 \textbf{Theorem \thesection.3} Let $A$ be a regular multiplier Hopf algebra
 with $M_{r}(V\otimes A) = M_{0}(V\otimes A)$ for any vector space $V$.
 Then the Yetter-drinfel'd module category ${}_{A}\mathcal{YD}^{A}$
 is equivalent to the centre $Z({}_{A}\mathcal{M})$ of left unital $A$-module category ${}_{A}\mathcal{M}$.

 \emph{Proof.}
 From above, $F$ and $G$ preserve the braiding. Also $FG=1$ and $GF=1$.
 So we have established the equivalence between $Z({}_{A}\mathcal{M})$ and ${}_{A}\mathcal{YD}^{A}$.
 $\hfill \blacksquare$
 \\

 The assumption $M_{r}(V\otimes A) = M_{0}(V\otimes A)$ as vector spaces in Proposition \thesection.2
 makes $C'_{A, V}(1\otimes v)$ belong to $M_{0}(V\otimes A)$, which is exactly what we need.
 This assumption sometimes naturally holds.

 If $A$ has a unit, $M_{r}(V\otimes A) = V\otimes A = M_{0}(V\otimes A)$, the assumption holds,
 then Theorem \thesection.3 actually shows that for a Hopf algebra $A$
 the Yetter-drinfel'd module category ${}_{A}\mathcal{YD}^{A}$
 is equivalent to the centre $Z({}_{A}\mathcal{M})$ of $A$-module category ${}_{A}\mathcal{M}$.

 If $A$ is commutative, the result $M_{r}(V\otimes A) = M_{0}(V\otimes A)$ as vector spaces also holds,
 and we get the following corollary

 \textbf{Corollary \thesection.4} Let $A$ be a commutative regular multiplier Hopf algebra.
 Then the Yetter-drinfel'd module category ${}_{A}\mathcal{YD}^{A}$
 is equivalent to the centre $Z({}_{A}\mathcal{M})$ of ${}_{A}\mathcal{M}$.
 \\

 From the proof in Proposition \thesection.2, $C_{X, V}$ is left and right $A$-linear, and these two structure are compatible.
 We can define $a\cdot C'_{A, V}(1\otimes v)$ by
 \begin{eqnarray*}
 (a\cdot C'_{A, V}(1\otimes v)) \cdot b
 &=& a\cdot (C'_{A, V}(1\otimes v) \cdot b)\\
 &=& a\cdot C_{A, V}(b\otimes v)
  =  C_{A, V}(a_{(1)}b\otimes a_{(2)}\cdot v) \\
 &=& C_{A, V}(a_{(1)}\otimes a_{(2)}\cdot v) \cdot b.
 \end{eqnarray*}
 Then $a\cdot C'_{A, V}(1\otimes v) = C_{A, V}(a_{(1)}\otimes a_{(2)}\cdot v) \in V\otimes A$.
 If $V$ is finite-dimensional, $C'_{A, V}(1\otimes v) \in M_{r}(V\otimes A)=V\otimes L(A)$,
 then  $C'_{A, V}(1\otimes v) \in V\otimes M(A) = M_{0}(V\otimes A)$.
  
 Denote by ${}_{A}\mathcal{YD}_{f.d}^{A}$ (${}_{A}\mathcal{M}_{f.d}$) the subcategory of ${}_{A}\mathcal{YD}^{A}$
 (${}_{A}\mathcal{M}$), in which the objects are all finite-dimensional.
 Then we have

 \textbf{Theorem \thesection.5} Let $A$ be a regular multiplier Hopf algebra.
 Then the category ${}_{A}\mathcal{YD}_{f.d}^{A}$
 is equivalent to the centre $Z({}_{A}\mathcal{M}_{f.d})$ of the category ${}_{A}\mathcal{M}_{f.d}$.

\section{Generalized Yetter-Drinfel'd modules over a regular multiplier Hopf algebra}
\def\theequation{\thesection.\arabic{equation}}
\setcounter{equation}{0}

 Using the extended modules (of the right modules) shown in section 2,
 we first give the generalized definition of (left-right) Yetter-Drinfel'd module
 over a regular multiplier Hopf algebra, which generalizes the notion before (e.g. in \cite{YW11}).
 We call it generalized Yetter-Drinfel'd module, and also we can define it by the extended modules of the right ones similarly.
 \\

 \textbf{Definition \thesection.1}
 Let $(A, \Delta, \varepsilon, S)$ be a regular multiplier Hopf algebra
 and $V$ a vector space. Then $V$ is called a generalized (left-right) Yetter-Drinfel'd module over $A$, if
 the following conditions hold:
 \begin{itemize}

 \item[(1)] $(V, \cdot)$ is a unital left $A$-module, i.e., $A\cdot V = V$;

 \item[(2)] $(V, \Gamma)$ is a right $A$-comodule, where $\Gamma: V\rightarrow M_{r}(V \otimes A)$ denotes the right coaction of $A$ on $V$
 and $M_{r}(V \otimes A)$ is the extended module of right $A$-module $V\otimes A$;

 \item[(3)] $\Gamma$ and $\cdot$ satisfy the following compatible condition:
 \begin{eqnarray}
 && (a\cdot v)_{(0)} \otimes (a\cdot v)_{(1)} a' = a_{(2)}\cdot v_{(0)} \otimes a_{(3)} v_{(1)} S^{-1}(a_{(1)}) a'. \label{1}
 \end{eqnarray}
 for all $a, a'\in A$ and $v\in V$.
 \end{itemize}

 \textbf{Remark \thesection.2}
 (1) Although the compatible conditions are same,
 the Yertter-Drinfel'd module in \cite{YW11} is different from this one, and can be considered as a special case of the above definition.
 Indeed, $M_{0}(V \otimes A)$, the range of coaction,
 is the intersection of extended modules of right and left $A$-module $V \otimes A$ (see section 2.3 in \cite{V08} for more details).

 (2) The compatible condition make sense, since
 $(a\cdot v)_{(0)} \otimes (a\cdot v)_{(1)} a' = \Gamma(a\cdot v)(1 \otimes a') \in V \otimes A$,
 and the right-hand side also belongs to $V \otimes A$. Indeed,
 $S^{-1}(a_{(1)}) a' \otimes a_{(2)} \otimes a_{(3)} \in A \otimes M(A \otimes A)$, and
 $v_{(0)} \otimes v_{(1)} S^{-1}(a_{(1)}) a' = \Gamma(v)(1 \otimes S^{-1}(a_{(1)}) a') \in V \otimes A$,
 then $a_{(2)}\cdot v_{(0)} \otimes a_{(3)} v_{(1)} S^{-1}(a_{(1)}) a' \in V \otimes A$.

 (3) The compatible condition (\ref{1}) is equivalent to the following condition
 \begin{eqnarray*}
 (a_{(2)}\cdot v)_{(0)} \otimes (a_{(2)}\cdot v)_{(1)} a_{(1)}a' = a_{(1)}\cdot v_{(0)} \otimes a_{(2)} v_{(1)}a'.
 \end{eqnarray*}

 (4) Let $A$ be a Hopf algebra with unit $1$, then $M_{r}(V \otimes A) = V\otimes A$. In this situation,
 the Yetter-Drinfel'd module for a Hopf algebra is a special case of
 Yetter-Drinfel'd module for a multiplier Hopf algebra. We will consider this at the end of this section.
 \\

 By the definition of generalized Yetter-Drinfel'd modules, we can define generalized Yetter-Drinfel'd module categories ${}_{A}\mathcal {GYD}^{A}$.
 The objects in ${}_{A}\mathcal{GYD}^{A}$ are generalized left-right Yetter-Drinfel'd modules,
 and the morphisms are linear maps which interwine with the
 left action and the right coaction of $A$ on $V$, i.e., the morphisms between two objects are
 left $A$-linear and right $A$-colinear maps. More precisely,
 let $V, W \in {}_{A}\mathcal{GYD}^{A}$ and $f: V\rightarrow W$ be a morphism, then
 \begin{eqnarray*}
 && f(a\cdot v) = a\cdot f(v), \\
 && (\Gamma_{W}\circ f) (v) (1\otimes a') = (f\otimes \iota)\Gamma_{V}(v)(1\otimes a'),
 \end{eqnarray*}
 for all $a, a' \in A$ and $v \in V$, where $\Gamma_{W}$ ($\Gamma_{V}$) is the right coaction on $W$ ($V$).

 The other three kind of generalized Yetter-Drinfel'd module categories can be also defined as in \cite{YW11}.
 \\

 The results in \cite{YW11} can be naturally generalized to this case, such as the following.

 \textbf{Theorem \thesection.3} The category ${}_{A}\mathcal {GYD}^{A}$ is a braided monoidal category.
 \\

 Denote by $Aut(A)$ the set of all isomorphisms $\alpha$ from $A$ to $A$ that are algebra maps satisfying
 $(\Delta \circ \alpha)(a)=(\alpha \otimes \alpha)\circ \Delta(a)$ for all $a\in A$.
 Then for $\alpha, \beta \in Aut(A)$, we also can define generalized $(\alpha, \beta)$-Yetter-Drinfel'd modules by the extended modules,
 and construct a new class of braided $T$-categories. The procedure is totally same as in \cite{V08}.

 Let $(A, \Delta, \varepsilon, S)$ be a regular multiplier Hopf algebra and $V$ a vector space.
 Then $V$ is called a (left-right) generalized $(\alpha, \beta)$-Yetter-Drinfel'd module over $A$, if
 \begin{itemize}

 \item[(1)]
  $(V, \cdot)$ is an unital left $A$-module, i.e., $A\cdot V = V$.

 \item[(2)]
  $(V, \Gamma)$ is a (right) $A$-comodule, where $\Gamma: M\rightarrow M_{r}(V \otimes A)$ denotes the right coaction of $A$ on $V$,
 and $M_{r}(V \otimes A)$ denote the extended module of right $A$-module $V \otimes A$.

 \item[(3)]
  $\Gamma$ and $\cdot$ satisfy the following compatible conditions
 \begin{eqnarray}
 && (a\cdot v)_{(0)} \otimes (a\cdot v)_{(1)} a' = a_{(2)}\cdot v_{(0)} \otimes \beta(a_{(3)}) v_{(1)} \alpha (S^{-1}(a_{(1)})) a'.
 \end{eqnarray}
 \end{itemize}

 Similarly, we can define the generalized $(\alpha, \beta)$-Yetter-Drinfel'd module category,
 denote it as $_{A}\mathcal{GYD}^{A}(\alpha, \beta)$,
 by which we can also construct a new braided $T$-category.
 Define $\mathcal{GYD}(A)$ as the disjoint union of all ${}_{A}\mathcal{GYD}^{A}(\alpha, \beta)$
 with $(\alpha, \beta) \in G$, where $G= Aut(A)\times Aut(A)$ is the twisted semi-direct square of group $Aut(A)$
 with product $ (\alpha, \beta)\# (\gamma, \delta)=(\alpha\gamma, \delta\gamma^{-1}\beta\gamma)$.
 Then we have
 \\

 \textbf{Theorem \thesection.4} $\mathcal{GYD}(A)$ is a braided $T$-category with the structures as follows.
 \begin{itemize}
 \item Tensor product:
 if $V \in {}_{A}\mathcal{GYD}^{A}(\alpha, \beta)$ and $W \in {}_{A}\mathcal{GYD}^{A}(\gamma, \delta)$
 with $\alpha, \beta , \gamma, \delta\in Aut(A)$,
 then $V\otimes W \in {}_{A}\mathcal{GYD}^{A} (\alpha\gamma, \delta\gamma^{-1}\beta\gamma) $, with the
 structures as follows:
 \begin{eqnarray*}
 && a\cdot (v\otimes  w) = \gamma(a_{(1)})\cdot v \otimes  \gamma^{-1}\beta\gamma(a_{(2)})\cdot w, \\
 && \Gamma(v\otimes  w)(1\otimes  a') = (v\otimes  w)_{(0)} \otimes  (v\otimes  w)_{(1)} a' = (v_{(0)} \otimes  w_{(0)}) \otimes  w_{(1)}v_{(1)} a'
 \end{eqnarray*}
 for all $v\in V, w\in W$ and $a, a'\in A$.

 \item Crossed functor:
 Let $W\in {}_{A}\mathcal{GYD}^{A}(\gamma, \delta)$, define $\varphi_{(\alpha, \beta)}(W)={}^{(\alpha, \beta)}W= W$
 as vector space, with structures: for all $a, a'\in A$ and $w\in W$
 \begin{eqnarray*}
 && a\rightharpoonup w = \gamma^{-1} \beta\gamma\alpha^{-1} (a)\cdot w, \\
 && \Gamma(w)(1\otimes  a')=: w_{<0>} \otimes  w_{<1>} a'= w_{(0)}\otimes \alpha\beta^{-1}(w_{(1)}) a'.
 \end{eqnarray*}
 Then
 ${}^{(\alpha, \beta)}W \in {}_{A}\mathcal{GYD}^{A}((\alpha, \beta)\# (\gamma, \delta)\# (\alpha, \beta)^{-1})
 = {}_{A}\mathcal{GYD}^{A}(\alpha\gamma\alpha^{-1}, \alpha\beta^{-1}\delta\gamma^{-1}\beta\gamma\alpha^{-1})$.
 The functor $\varphi_{(\alpha, \beta)}$ acts as identity on morphisms.

 \item Braiding:
 If $V \in {}_{A}\mathcal{GYD}^{A}(\alpha, \beta)$, and $W \in {}_{A}\mathcal{GYD}^{A}(\gamma, \delta)$.
 Take ${}^{V}W={}^{(\alpha, \beta)}W$, define a map $C_{V, W}: V\otimes W \longrightarrow {}^{V}W \otimes V$ by
 \begin{eqnarray*}
 && C_{V, W}(v\otimes w) = w_{(0)} \otimes \beta^{-1}(w_{(1)})\cdot v
 \end{eqnarray*}
 for all $v\in V$ and $w\in W$.
 \end{itemize}
 \bigskip

 When $A$ has a unit, i.e., $A$ is a Hopf algebra, $M_{r}(V \otimes A) = V \otimes A$,
 and the $_{A}\mathcal {GYD}^{A}(\alpha, \beta)$ is exactly $_{A}\mathcal {YD}^{A}(\alpha, \beta)$ in the Hopf algebra case.
 Especially when $A$ is a coFrobenius Hopf algebra with cointegral $t$ in $A$ and integral $\varphi$ on $A$ satisfying $\varphi(t)=1$.
 Then we have
 \begin{eqnarray}
 {}_{A}\mathcal{YD}^{A}(\alpha, \beta) \cong {}_{\widehat{A} \bowtie A(\alpha, \beta)}\mathcal{M}, \label{3}
 \end{eqnarray}
 where ${}_{\widehat{A} \bowtie A(\alpha, \beta)}\mathcal{M}$ is the unital $\widehat{A} \bowtie A(\alpha, \beta)$-module category,
 and $\widehat{A} \bowtie A (\alpha, \beta)$ is the diagonal crossed product with the multiplication (see Definition 2.4 in \cite{YW11a})
 $$(p\bowtie a)(p'\bowtie a') =  p(\alpha(a_{(1)})\blacktriangleright p' \blacktriangleleft S^{-1}\beta(a_{(3)}))\otimes a_{(2)}a'$$
 for $a, a'\in A$ and $p, p' \in \widehat{A}$.

 Indeed, the correspondence is give as follows.
 If $M \in {}_{A}\mathcal{YD}^{A}(\alpha, \beta)$, then $M \in {}_{\widehat{A} \bowtie A(\alpha, \beta)}\mathcal{M}$
 with the structure
 $$(p\bowtie a) \cdot m = p((a \cdot m)_{(1)}) (a \cdot m)_{(0)}.$$
 Conversely, if $M \in {}_{\widehat{A} \bowtie A(\alpha, \beta)}\mathcal{M}$, then $M \in {}_{A}\mathcal{YD}^{A}(\alpha, \beta)$
 with structures
 \begin{eqnarray*}
 && a\cdot m = (\varepsilon \bowtie a) \cdot m, \\
 && m\mapsto m_{(0)} \otimes m_{(1)} = (\varphi(\cdot t_{(2)}) \bowtie 1)\cdot m \otimes S^{-1}(t_{(1)}).
 \end{eqnarray*}

 If $\alpha=\beta=\iota$, the result (\ref{3}) is just Theorem 12 shown in \cite{Zh99}.

\section{Yetter-Drinfel'd $(H, A)$-modules over a regular multiplier Hopf algebra}
\def\theequation{\thesection.\arabic{equation}}
\setcounter{equation}{0}

 Let $A$ be a regular multiplier Hopf algebra, and $R$ be an algebra with or without identity, but with non-degenerate product.
 Assume that $R$ is an unital left $A$-module, we recall from \cite{DrVZ99,VZ99} that $R$ is called a (left) $A$-module algebra if
 \begin{eqnarray}
 a \cdot (xx') = (a_{(1)} \cdot x) (a_{(2)} \cdot x') \label{2}
 \end{eqnarray}
 for all $a\in A$ and $x, x' \in R$.
 We can extend the action of $A$ on $R$ to $M(R)$ by the following formulas
 \begin{eqnarray}
 && (a \cdot x)x' = a_{(1)} \big(x (S(a_{(2)}) \cdot x') \big), \\
 && x(a \cdot x') = a_{(2)} \big((S^{-1}(a_{(1)}) \cdot x) x' \big).
 \end{eqnarray}
 It follows that $M(R)$ is a left $A$-module, the action of $A$ on $M(R)$ is still non-degenerate but no longer unital,
 and $a \cdot 1 = \varepsilon(a) 1$.
 \\

 By a right coaction of $A$ on $R$, we mean an injective map $\Gamma: R\longrightarrow M(R \otimes A)$ satisfying
 \begin{itemize}
   \item[(i)] $\Gamma(R)(1\otimes A)\subseteq R\otimes A$ and $(1\otimes A)\Gamma(R)\subseteq R\otimes A$,
   \item[(ii)] $(\Gamma\otimes\iota)\Gamma=(\iota\otimes\Delta)\Gamma$.
 \end{itemize}
 If furthermore $\Gamma$ is a homomorphism, then $R$ is called a (right) $A$-comodule algebra.
 The injectivity of $\Gamma$ is equivalent to the counitary property $(\iota \otimes \varepsilon)\Gamma = \iota$.
 \\

 A (left-right) Yetter-Drinfel'd $A$-module algebra $R$ is an algebra $R$ (with or without unit) over the field $K$,
 which is an unital left $A$-module algebra and a right $A$-comodule algebra satisfying the compatibility condition (\ref{1}).
 In the another word, $R$ is an algebra (with or without the unitary property) in Yetter-Drinfel'd module category $\mathcal{Q}^{A}$.

 A Yetter-Drinfel'd $A$-module algebra $R$ is called \emph{$A$-commutative},
 if for all $x, y\in R$,
 \begin{eqnarray*}
 xy = y_{(0)}(y_{(1)} \cdot x).
 \end{eqnarray*}

 Let $A$ be a quasitriangular multiplier Hopf algebra with quasitriangular structure $\mathcal{R}$  and $H$ be a left $A$-module algebra.
 Then there is a natural right $A$-comodule structure on $H$, namely
 \begin{equation}\label{3}
 \rho(h) = \tau(\mathcal{R})(h\otimes 1).
 \end{equation}
 Remark here that (\ref{3}) makes sense. Since the flip map $\tau: A\otimes A \longrightarrow A\otimes A $ is a homomorphism,
 it can be extended to $M(A\otimes A)$. For $h\in H$, there exists finite $a_{i}$ and $h_{i}$ such that $h=\sum_{i}a_{i}\cdot h_{i}$.
 So $\tau(\mathcal{R})(h\otimes 1)=\sum_{i}\tau(\mathcal{R})(a_{i}\otimes 1)(h_{i}\otimes 1)\in M(R\otimes A)$.
 \\

 \textbf{Proposition \thesection.1}
 Let $A$ be a quasitriangular multiplier Hopf algebra and $H$ be a left $A$-module algebra with a $A$-comodule structure as (\ref{3}).
 Then $H$ is a right $A$-comodule algebra, and a Yetter-Drinfel'd $A$-module algebra.
 \\

 As shown in \cite{YW11}, the (left-right) Yetter-Drinfel'd module category $\mathcal{Q}^{A}$ is a braided monoidal category.
 From the above proposition, if $A$ is quasitriangular,
 then the unital $A$-module category ${}_{A}\mathcal{M}$ is a braided monoidal subcategory of ${}_{A}\mathcal {YD}^{A}$,
 and the braiding on ${}_{A}\mathcal{M}$ is giveb by
 $C_{M, N}: M\otimes N \longrightarrow N\otimes M$, $C_{M, N}(m\otimes n)=\tau(\mathcal{R})(n\otimes m)$.
 \\

 Let $A$ be a regular multiplier Hopf algebra and $H$ be a Yetter-Drinfel'd module algebra.
 In the following, we introduced a new kind of modules: Yetter-Drinfel'd $(H, A)$-modules.
 \\

 \textbf{Definition \thesection.2}
 Let $H$ be a Yetter-Drinfel'd $A$-module algebra. A Yetter-Drinfel'd $(H, A)$-module $M$ is a $K$-module, which has the structures
 of Yetter-Drinfel'd ($A$-)module and of left $H$-module with left $H$-action (denoted by $\rightharpoonup$) such that
 the following compatibility conditions hold: for all $h\in H$, $m\in M$ and $a, a'\in A$,
 \begin{eqnarray}
 && a\cdot (h\rightharpoonup m) = (a_{(1)} \cdot h) \rightharpoonup (a_{(2)} \cdot m), \\
 && \rho(h\rightharpoonup m)(1\otimes a') = h_{(0)}\rightharpoonup m_{(0)} \otimes m_{(1)}h_{(1)}a'.
 \end{eqnarray}
 Denote by ${}_{H}\mathcal{Q}^{A}$ the category of Yetter-Drinfel'd $(H, A)$-modules and Yetter-Drinfel'd $(H, A)$-module maps.
 \\

 Note that if $H$ is the trivial Yetter-Drinfel'd module algebra $K$, then the Yetter-Drinfel'd $(H, A)$-module is nothing but
 Yetter-Drinfel'd $A$-module.
 If $A$ is a coFrobenius Hopf algebra, by \cite{Zh99} $\mathcal{Q}^{A}\cong {}_{D(A)}\mathcal{M}$.
 If $H$ ia a $D(A)$-module algebra, then ${}_{H}\mathcal{Q}^{A}$ is equivalent to the unital module category
 ${}_{H\# D(A)}\mathcal{M}$, where $D(A)=\widehat{A}\bowtie A$ the Drinfel'd double of $A$, and $H\# D(A)$ the usual smash product.
 \\

 Now we consider the case in which $H$ ia an unital $A$-commutative Yetter-Drinfel'd module algebra. Take $M$ in ${}_{H}\mathcal{Q}^{A}$,
 we may define a right action of $H$ on $M$ as follows: for any $m\in M$ and $h\in H$,
 \begin{equation}\label{4}
   m \leftharpoonup h = h_{(0)} \rightharpoonup (h_{(1)} \cdot m)
 \end{equation}
 It is no hard to check that this action is well-defined and along with the left $H$-module structure of $M$
 makes $M$ into an $H$-bimodule.

 For $M, N \in {}_{H}\mathcal{Q}^{A}$, we form a tensor product $M\otimes_{H} N$, and endow it with the following $H$-action
 and $A$-structures: for any $m\in M$, $n\in N$, $a, a' \in A$ and $h\in H$,
 \begin{eqnarray*}
 && a \cdot (m\otimes n) = a_{(1)} \cdot m \otimes a_{(2)} \cdot n, \\
 && \chi(m\otimes n)(1 \otimes a') = m_{(0)} \otimes n_{(0)} \otimes n_{(1)}m_{(1)} a', \\
 && h \rightharpoonup (m\otimes n) = h\rightharpoonup m \otimes n.
 \end{eqnarray*}
 Then $M\otimes_{H} N$ with the above structures is a Yetter-Drinfel'd $(H, A)$-module, denoted by $M\widetilde{\otimes}_{H} N$.

 Note that the right $H$-module structure of $M\widetilde{\otimes}_{H} N$, defined as in (\ref{4}),
 is the same as the one coming from the right $H$-module $N$, i.e. $(m\otimes n)\leftharpoonup h = m\otimes n\leftharpoonup h$
 for $m\in M$, $n\in N$ and $h\in H$. Therefore, the standard map for Yetter-Drinfel'd $(H, A)$-modules $X$, $Y$ and $Z$,
 \begin{equation*}
 \Phi: (X\widetilde{\otimes}_{H} Y)\widetilde{\otimes}_{H} Z \longrightarrow  X\widetilde{\otimes}_{H} (Y \widetilde{\otimes}_{H} Z)
 \end{equation*}
 is $H$-linear, and hence an isomorphism in ${}_{H}\mathcal{Q}^{A}$.
 By the construction above, we have

 \textbf{Theorem \thesection.3}
 With notation above, $({}_{H}\mathcal{Q}^{A}, \widetilde{\otimes}_{H}, H)$ forms a monoidal category.

\section*{Acknowledgements}

 The work was partially supported by the NNSF of China (No. 11226070), the NJAUF (No. LXY201201019, LXYQ201201103).
 The first author Tao Yang is very grateful to Professor A. Van Daele for his lectures and discussions on the multiplier Hopf algebra theory.

\end {document}